\documentclass[a4paper,12pt]{amsart}
\usepackage[width=0.8\paperwidth,height=0.85\paperheight,twoside=false,includehead,includefoot]{geometry}
\usepackage{amsmath,amssymb,amsfonts}
\usepackage{amsthm,graphicx}


\theoremstyle{plain}
\newtheorem{thm}{Theorem}[section]

\newtheorem{conj}{Conjecture}
\newtheorem{ex}[thm]{Example}
\newtheorem{rem}[thm]{Remark}

\def\N{\mathbb {Z}_{\geq 1}} \def\Z{\mathbb {Z}} \def\Q{\mathbb {Q}}
\def\R{\mathbb {R}} 
\def\d{\partial} 
   
\def\hh{\ast} 
\def\A{\mathcal{A}}
 
 \def\l{\textrm{\textup{\textmd{\textbf{l}}}}} \def\m{\textrm{\textup{\textmd{\textbf{m}}}}}

\def\z{\zeta}

\def\ZZ{\mathcal {Z}}
\def\fh{\frak H}

\def\za{\zeta_{\mathcal{A}}} 
\def\zff{\zeta_{\mathcal{F}}}
\def\zf{\zeta_{\mathcal{S}}} 

\font\fivecy=wncyr5 \def\sa{\hbox{\fivecy X}}
 
\def\ZA{\mathit{Z}_{\mathcal{A}}}
\def\ZS{\mathit{Z}_{\mathcal{S}}}
\def\ZF{\mathit{Z}_{\mathcal{F}}}

\begin{document}
\title[Derivation relations for FMZVs]{Derivation relations for finite multiple zeta values}
\author{Hideki Murahara}
\date{2016.2.26}
\address{Graduate School of Mathematics, Kyushu University \\
744 Motooka Nishi-ku, Fukuoka, 819-0395 Japan,}
\email{h-murahara@gmail.com} 

\begin{abstract}
The derivation relations for multiple zeta values is proved by Ihara, Kaneko and Zagier. 
We prove its counterpart for finite multiple zeta values.   
\end{abstract}
\maketitle

\section{Introduction}
For integers $k_1, \ldots , k_{d}\in \N$ with $k_1 \geq 2$, the multiple zeta value (MZV) is defined by 
\begin{align*}
\z (k_1,\ldots, k_d) := \sum_{n_1>\cdots >n_d \geq 1} \frac {1}{n_1^{k_1}\cdots n_d^{k_d}}.
\end{align*} 
To describe the derivation relations for MZVs, we use the algebraic setup introduced by M.\ Hoffman \cite{hoffman-97}. 
%%%%%
Let $ \fh =\Q \left\langle x,y \right\rangle$ be the noncommutative polynomial ring in two indeterminates $x$, $y$ and $\fh^1$ and $\fh^0$ its subrings $\Q +\fh y$ and $\Q +x \fh y$. 
We set $z_{k} = x^{k-1} y$ $(k\in\N)$. 
Then $\fh^1$ is freely generated by $\{z_{k}\}_{k \geq 1}$. 
For any word $w$, let $|w|$ be the total degree.

We define the $\Q$-linear map $\mathit{Z}:\fh^0 \to \R$ by $\mathit{Z}(1):=1$ and $\mathit{Z}( z_{k_1} \cdots z_{k_d}):= \z (k_1,\ldots, k_d)$, and the harmonic product $\hh$ on $\fh^1$ inductively by
\begin{align*}
&\quad\quad\,\,\,\,\,\,  1\hh w= w\hh 1 = w , \\
z_{k} w_{1} \hh z_{l} w_{2}&= z_{k} (w_{1} \hh z_{l} w_{2}) + z_{l} (z_{k} w_{1} \hh w_{2}) + z_{k+l} (w_{1} \hh w_{2}) 
\end{align*}
($k,l \in \mathbb{Z}_{\geq 1}$ and $w$, $w_{1}$, $w_{2}$ are words in $\fh^1$), together with $\mathbb{Q}$-bilinearity. 
The harmonic product $\hh$ is commutative and associative, therefore $\fh^1$ is $\Q$-commutative algebra with respect to $\hh$. 
%We denote it by $\mathfrak{H}^{1}_{\hh}$. 
The subset $\fh^0$ is a subalgebra of $\fh^1$ with respect to $\hh$.
% and we denote it by $\fh^{0}_{\hh}$. 
%We then have 
%\[Z(w_{1} \hh w_{2}) = Z(w_{1}) Z(w_{2}) \]
%for any $w_{1}, w_{2} \in \fh^0$. 

A derivation $\d$ on $\fh$ is a $\Q$-linear endomorphism of $\fh$ satisfying Leibniz's rule $\partial(ww^{\prime})=\partial(w)w^{\prime}+w\partial(w^{\prime})$. Such a derivation is uniquely determined by its images of generators $x$ and $y$. Set $z:=x+y$. 
For each $l\ge 1$, the derivation $\d_l:\fh \to \fh$ is defined by $\partial_l(x):=xz^{l-1}y$ and $\d_l(y):=-xz^{l-1}y$. 
We note that $\d_l(1)=0$ and $\d_l(z)=0$. 
In \cite{ihara_kaneko_zagier-2006}, K.\ Ihara, M.\ Kaneko and D.\ Zagier proved the derivation relations for MZVs. 
\begin{thm}[Ihara--Kaneko--Zagier] \label{ihara-kaneko-zagier}
For $l\in\N$, we have 
\begin{align*} 
 \mathit{Z} (\partial_{l}(w)) = 0 \quad (w \in \fh^0).  
\end{align*}
\end{thm}
In this paper, we prove its counterpart for what we call `finite multiple zeta values', a generic term for $\A$-finite multiple zeta values and symmetrized multiple zeta values, which we now explain.    

We consider the collection of truncated sums $\z_{p}(k_1, \ldots, k_d):=\sum_{p>n_1>\cdots>n_d\geq1}\frac{1}{n_1^{k_1}\cdots n_d^{k_d}}$ modulo all primes $p$ in the quotient ring $\A=(\prod_{p}\Z/p\Z)/(\bigoplus_{p}\Z/p\Z)$, which is a $\Q$-algebra. 
Elements of $\A$ are represented by $(a_p)_p$, where $a_p\in\Z/p\Z$, and two elements $(a_p)_p$ and $(b_p)_p$ are identified if and only if $a_p=b_p$ for all but finitely many primes $p$.
For integers $k_1, \ldots ,k_{d} \in \N$, the $\A$-finite multiple zeta value ($\A$-FMZV) $\za(k_1,\cdots, k_d)$ is defined by
\begin{align*}
 \za (k_1, \ldots, k_d)&:=\biggl( \sum_{p>n_1>\cdots >n_d\geq1}\frac{1}{n_1^{k_1}\cdots n_d^{k_d}} \bmod{p} \biggr)_{p} \in \A.
\end{align*}
We denote by $\ZZ_{\A}$ the $\Q$-vector subspace of $\A$ spanned by $1$ and all $\A$-finite multiple zeta values. 
It is known that this is also a $\Q$-algebra.

The symmetrized multiple zeta values or finite real multiple zeta values, which were first introduced by Kaneko--Zagier \cite{kaneko_zagier-2014, kaneko_zagier-2015}, are defined for any positive integers $k_1,\ldots ,k_{d}$ as follows:  
\begin{align*}
\zf ^{\ast} (k_1, \ldots, k_d) &:= \sum_{i=0}^{d}(-1)^{k_1+\cdots +k_i}\z^{\ast}(k_i, \ldots , k_1) \z ^{\ast}(k_{i+1}, \ldots , k_d), \\
\zf ^{\sa} (k_1, \ldots, k_d) &:= \sum_{i=0}^{d}(-1)^{k_1+\cdots +k_i}\z^{\sa}(k_i, \ldots , k_1) \z ^{\sa}(k_{i+1}, \ldots , k_d). 
\end{align*}
Here, the symbols $\z ^{\ast}$ and $\z ^{\sa}$ on the right-hand sides stand for the regularized values coming from harmonic and shuffle regularizations respectively,
i.e.,  real values obtained by taking constant terms of harmonic and shuffle regularizations as explained in   \cite{ihara_kaneko_zagier-2006}. 
In the sums, we understand $\zeta^{\ast}(\emptyset )=\zeta^{\sa}(\emptyset)=1$.

Let $\ZZ_{\R}$ be the $\Q$-vector subspace of $\R$ spanned by $1$ and all MZVs. 
It is known that this is a $\Q$-algebra.
In \cite{kaneko_zagier-2014, kaneko_zagier-2015}, Kaneko and Zagier proved that the difference $\zf^{\ast} (k_1,\ldots , k_d) -\zf^{\sa} (k_1, \ldots ,k_d)$
is in the principal ideal of $\ZZ_{\R}$ generated by $\zeta(2)$ (or $\pi^2$), in other words, that the congruence   
\[  \zf^{\ast} (k_1, \ldots , k_d) \equiv \zf^{\sa} (k_1, \ldots , k_d)  \pmod{\zeta (2)} \]
holds in $\ZZ_{\R}$. 
They then defined the symmetrized multiple zeta value (SMZV)  $\zf (k_1, \ldots , k_d)$ as an element 
in the quotient ring $\ZZ_{\R}/\zeta (2)$ by 
\[  \zf (k_1, \ldots , k_d) := \zf^{\ast} (k_1, \ldots , k_d) \bmod \zeta(2). \]
We also refer to the values $\zf^{\ast} (k_1,\ldots ,k_d)$ and $\zf^{\sa} (k_1\ldots , k_d)$ as
(harmonic and shuffle versions of) symmetrized multiple zeta values.

Then, Kaneko and Zagier conjectured the following:  
\begin{conj}[Kaneko--Zagier]
There exists an algebra isomorphism $\phi$ between $\ZZ_{\A}$ and $\ZZ_{\R}/\z(2)$ such that 
\begin{equation*}
\phi :
\begin{matrix}
\ZZ_{\A} & \rightarrow & \ZZ_{\R}/\z(2) \\
\text{\rotatebox{90} {$\in$}}&    &\text{\rotatebox{90} {$\in$}} \\
\za (k_1,\ldots, k_d) & \mapsto & \zeta_{\mathcal{S}} (k_1, \ldots, k_d).
\end{matrix}
\end{equation*}
\end{conj}
We define two $\Q$-linear maps $\ZA \colon\fh^1 \to \A$ and $\ZS \colon\fh^1 \to \ZZ_{\R}/\zeta (2)$ by $\ZA (1):=1$ and $\ZA (z_{k_1} \cdots z_{k_d}):=\za (k_1,\ldots ,k_d)$, and $\ZS (1):=1$ and $\ZS (z_{k_1} \cdots z_{k_d}):=\zeta_{\mathcal{S}} (k_1,\ldots ,k_d)$, respectively. 

We finish this section by mentioning the harmonic product rule and the duality theorem for FMZVs. 
The former is due to Hoffman \cite{hoffman-97} for $\A$-FMZVs, and Kaneko and Zagier \cite{kaneko_zagier-2014,kaneko_zagier-2015} for SMZVs. 
We use these results in the proof of our main theorem. 
\begin{thm}[Hoffman, Kaneko--Zagier]  \label{3.3}
For any words $w, w'\in\fh^1$, we have
\begin{align*}
\ZF (w\hh w')&=\ZF (w) \ZF (w'), 
\end{align*}
where the symbol `$\mathcal{F}$' stands either for `$\mathcal{A}$' or `$\mathcal{S}$'.
\end{thm}
The duality theorems for $\A$-finite and symmetrized versions are proved by Hoffman \cite{hoffman-2015} and D.\ Jarossay \cite{jarossay-2014}, respectively. 
\begin{thm} [Hoffman, Jarossay] \label{thm:duality}
Let $\phi$ be an automorphism on $\fh$ defined by 
\[ \phi (x)=z(=x+y),\quad \phi (y)=-y. \] 
Then, we have 
\[ \ZF (w)=\ZF (\phi (w)) \quad (w \in \fh^1). \]
\end{thm}

\section{Main theorems}
Kojiro Oyama conjectured the following derivation relations for FMZVs.
\begin{conj}[Oyama] \label{oyama}
For $l\in\N$, we have 
\begin{align} \label{1}
 \ZF (\partial_{l}(w)) = -\ZF (z^{l-1}yw) \quad (w \in \fh^1 ,\mathcal{F}=\mathcal{A} \textrm{ or } \mathcal{S}).  
\end{align}
\end{conj} 
Oyama proved this for $l\leq 4$ and Mitsuki Kosaki extended the proof further to $l\leq 6$. 
The aim of this paper is to prove this conjecture for all $l$. 
Actually, we prove the identity in the following form, which looks more general than the conjecture but in fact is equivalent to the conjecture. 
The proof of Theorem \ref{derivation0} will be given in the next section. 
\begin{thm} \label{derivation0}
For $\m =(m_1,\ldots ,m_s)\in (\N)^s \,\, (s\geq 0)$ and $l\in\N$, we have
\begin{align*}
&\ZF (z^{m_1-1}y \cdots z^{m_s-1}y \partial_{l}(w)) \\
&=-\ZF (z^{l-1}y z^{m_{1}-1}y \cdots z^{m_{s}-1}yw) \\
&\,\,\,\,\, +\sum_{i=1}^{s} \ZF (z^{m_1-1}y \cdots z^{m_{i-1}-1}y z^{m_i-1}x z^{l-1}y z^{m_{i+1}-1}y \cdots z^{m_{s}-1}yw) \quad (w \in \fh^1). 
\end{align*}
When $s=0$, we understand $z^{m_1-1}y \cdots z^{m_s-1}y=1$ on the left, and the right-hand side is $-\ZF (z^{l-1}yw)$, which yield Conjecture \ref{oyama}. 
\end{thm}
\begin{rem}
We see Conjecture \ref{oyama} implies Theorem \ref{derivation0} by putting $z^{m_1-1}y\cdots z^{m_s-1}yw$ for $w$ in eq.$(\ref{1})$, because 
\begin{align*}
\d_{l}(z^{m_1-1}y\cdots z^{m_s-1}yw) =&-z^{m_1-1}xz^{l-1}yz^{m_2-1}y\cdots z^{m_s-1}yw \\ 
&-\quad\cdots\cdots \\
&-z^{m_1-1}y\cdots z^{m_s-1}xz^{l-1}yw \\ 
&+z^{m_1-1}y\cdots z^{m_s-1}y\d_{l}(w)   
\end{align*}
by the definition of $\d_l$ \textrm{(}note $\d_l(z)=0$\textrm{)}, and 
\[ \ZF (\partial_{l}(z^{m_1-1}y \cdots z^{m_s-1}yw)) = -\ZF (z^{l-1}yz^{m_1-1}y \cdots z^{m_s-1}yw) \]
by eq.$(\ref{1})$. 
\end{rem}
\begin{ex}
When $l=3$ and $w=xy$ in Conjecture \ref{oyama}, we have 
\begin{align*}
\ZF (\d_3 (xy))&=-\ZF (z^2yxy) \\
&=-\ZF (x^2yxy+xy^2xy+yxyxy+y^3xy).
\end{align*}
Since
\begin{align*}
\d_3 (xy)&=xz^2y^2-x^2z^2y \\
&=xyxy^2+xy^4-x^4y-x^2yxy,
\end{align*}
we get 
\begin{align*}
\zff (5)-\zff (2,2,1)-\zff (2,1,2)-\zff (1,2,2)-\zff (2,1,1,1)-\zff (1,1,1,2)=0.    
\end{align*}
\end{ex}
\begin{ex}
The case $s=2$ in Theorem \ref{derivation0} gives 
\begin{align*}
\ZF (z^{m_1-1}yz^{m_2-1}y\d_l(w))&=-\ZF (z^{l-1}yz^{m_1-1}yz^{m_2-1}yw) \\
&\quad\,  +\ZF (z^{m_1-1}xz^{l-1}yz^{m_2-1}yw) \\ 
&\quad\,  +\ZF (z^{m_1-1}yz^{m_2-1}xz^{l-1}yw).  
\end{align*}
When $m_1=2, m_2=1, l=2$ and $w=y$, we get  
\begin{align*}
&\zff (4,1,1)+\zff (2,3,1)+\zff (2,1,3)+\zff (3,1,1,1)+\zff (1,3,1,1)+\zff (1,1,3,1)+\zff (1,1,1,3) \\
&+\zff (2,1,2,1)-\zff (2,1,1,1,1)+\zff (1,2,1,1,1)+\zff (1,1,1,2,1)-\zff (1,1,1,1,1,1)=0.    
\end{align*} 
\end{ex}
Let $S_n$ be the symmetric group of order $n$, which acts on any index  $\textrm{\textup{\textmd{\textbf{a}}}}=(a_1,\ldots ,a_n)$ by $\sigma (\textrm{\textup{\textmd{\textbf{a}}}})=(a_{\sigma (1)},\ldots ,a_{\sigma (n)})$. 
For an integer $s$ with $1\leq s\leq n$, let $S_n^{(s)}$ be the subset of $S_n$ given by
\[ S_{n}^{(s)}=\left\{ \sigma\in S_{n} \mid \sigma^{-1}(1)<\cdots <\sigma^{-1}(s) \right\} . \]
Under these notations, we have the following theorem, which is in fact an almost immediate consequence of Theorem \ref{derivation0}. 
The proof will also be given in the next section.  
\begin{thm} \label{derivation1}
For $\m =(m_1,\ldots ,m_s)\in (\N)^s \,\, (s\geq 0)$ and $\l =(l_1,\ldots ,l_t)\in (\N)^t \,\, (t\geq 1)$, 
we set $\textrm{\textup{\textmd{\textbf{a}}}}=(a_1,\ldots ,a_{s+t})=(\m ,\l)$. 
Then, we have  
\begin{align*} 
& \ZF (z^{m_1-1}y \cdots z^{m_s-1}y \partial_{l_1} \cdots \partial_{l_t}(w)) \\
& = (-1)^{s} \sum_{\sigma\in S_{s+t}^{(s)}} \ZF (z^{a_{\sigma (1)}-1}u^{\sigma}_1 \cdots z^{a_{\sigma (s+t)}-1} u^{\sigma}_{s+t}w) \quad (w \in \fh^1). 
\end{align*} 
Here, we set $u^{\sigma}_i=x$ if `$\sigma (i)\leq s$ and $\sigma (i+1)>s$' or `$\sigma (i)>s$ and $\sigma (i)<\sigma (i+1)$', and $u^{\sigma}_i=-y$ otherwise. 
\end{thm} 
\begin{ex}
When $s=1,t=2$ in Theorem \ref{derivation1}, we have 
\begin{align*}
\ZF (z^{m_1-1}y\d_{l_1}\d_{l_2}(w))&=\ZF (z^{m_1-1}xz^{l_1-1}xz^{l_2-1}yw)-\ZF (z^{m_1-1}xz^{l_2-1}yz^{l_1-1}yw) \\ 
&\,\,\,\,\, -\ZF (z^{l_1-1}yz^{m_1-1}xz^{l_2-1}yw)-\ZF (z^{l_2-1}yz^{m_1-1}xz^{l_1-1}yw) \\
&\,\,\,\,\, -\ZF (z^{l_1-1}xz^{l_2-1}yz^{m_1-1}yw)+\ZF (z^{l_2-1}yz^{l_1-1}yz^{m_1-1}yw). 
\end{align*}
By putting $ m_1=2,l_1=2, l_2=1$ and $w=y$, we get  
\begin{align*}
&\,\,\, \zff (5,1)-\zff (2,4)-\zff (3,2,1)-\zff (2,3,1)-\zff (1,1,4)\\ 
&-2\zff (3,1,1,1)-\zff (1,3,1,1)-2\zff (1,1,3,1)-\zff (2,1,2,1)+\zff (2,2,1,1) \\
&-\zff (1,2,1,1,1)+\zff (1,1,1,1,1,1)=0.  
\end{align*}
\end{ex}
\begin{rem}
For two indices $\m, \m'$, we say $\m'$ refines $\m$ (denoted  $\m' \succeq \m$) if $\m$ can be obtained from $\m'$ by combining some of its adjacent parts. 
%\preceq
Then, we have 
\begin{align}
\begin{split} \label{2}
&\ZF (x^{m_1-1}y \cdots x^{m_s-1}y \partial_{l_1} \cdots \partial_{l_t}(w)) \\
&= (-1)^{s} \sum_{\m' \succeq \m} \sum_{\sigma\in S_{s'+t}^{(s')}} \ZF (z^{a'_{\sigma (1)}-1}u^{\sigma}_1 \cdots z^{a'_{\sigma (s'+t)}-1} u^{\sigma}_{s'+t}w) \quad (w \in \fh^1), 
\end{split}  
\end{align}
where $\m ' =(m'_1,\ldots ,m'_{s'})$ and $\textrm{\textup{\textmd{\textbf{a}}}}'=(a'_1,\ldots ,a'_{s'+t})=(\m',\l)$. 
We note here that eq.$(\ref{2})$ is equivalent to Theorem \ref{derivation1}. 
Assume that Theorem \ref{derivation1} holds, we see by 
$x^{m_1-1}y\cdots x^{m_s-1}y=\sum_{\m' \succeq \m}  (-1)^{s'-s} z^{m'_1-1}y \cdots z^{m'_{s'}-1}y$ that 
\begin{align*}
&\ZF (x^{m_1-1}y\cdots x^{m_s-1}y\d_{l_1}\cdots\d_{l_t}(w)) \\
&=\sum_{\m' \succeq \m}  (-1)^{s'-s} \ZF (z^{m'_1-1}y \cdots z^{m'_{s'}-1}y\d_{l_1}\cdots\d_{l_t}(w)) \\
&= (-1)^{s} \sum_{\m' \succeq \m} \sum_{\sigma\in S_{s'+t}^{(s')}} \ZF (z^{a'_{\sigma (1)}-1}u^{\sigma}_1 \cdots z^{a'_{\sigma (s'+t)}-1} u^{\sigma}_{s'+t}w). 
\end{align*}
Conversely, assume that eq.$(\ref{2})$ holds.  
%Let $\m''=(m''_1,\ldots ,m''_{s''})$ and $\textrm{\textup{\textmd{\textbf{a}}}}''=(a''_1,\ldots ,a''_{s''+t})=(\m'',\l)$. 
Since 
\[ z^{m_1-1}y\cdots z^{m_s-1}y=\sum_{\m' \succeq \m} x^{m'_1-1}y \cdots x^{m'_{s'}-1}y, \]
and
\[ \sum_{\m' \succeq \m}  (-1)^{s'} \sum_{\substack{\m'' \succeq \m' \\ \m'' =(m''_1,\ldots ,m''_{s''})}} (m''_1,\ldots ,m''_{s''})=(-1)^{s} (m_1,\ldots ,m_s), \]
(the second equality is an identity of formal sums of indices) we have
\begin{align*}
&\ZF (z^{m_1-1}y\cdots z^{m_s-1}y\d_{l_1}\cdots\d_{l_t}(w)) \\
&=\sum_{\m' \succeq \m} \ZF (x^{m'_1-1}y \cdots x^{m'_{s'}-1}y\d_{l_1}\cdots\d_{l_t}(w)) \\
&=\sum_{\m' \succeq \m}  (-1)^{s'} \sum_{\m'' \succeq \m'} \sum_{\sigma\in S_{s''+t}^{(s'')}} \ZF (z^{a''_{\sigma (1)}-1}u^{\sigma}_1 \cdots z^{a''_{\sigma (s''+t)}-1} u^{\sigma}_{s''+t}w) \\
&=(-1)^{s} \sum_{\sigma\in S_{s+t}^{(s)}} \ZF (z^{a_{\sigma (1)}-1}u^{\sigma}_1 \cdots z^{a_{\sigma (s+t)}-1} u^{\sigma}_{s+t}w), 
\end{align*}
where $s''$ is the depth of $\m''$ and $\textrm{\textup{\textmd{\textbf{a}}}}'' =(a''_1,\ldots ,a''_{s''})=(\m'',\l)$.
\end{rem}

Before closing this section, we mention the maximal number of linearly independent relations supplied by Conjecture \ref{oyama}. 
In Table $1$, the first line means the weight of FMZVs 
(we call $k:=k_1+\cdots +k_d$ the weight for $\zff (k_1\ldots ,k_d)$). 
The second line gives the number of linearly independent elements in $\fh$ among all $\d_l(w)+z^{l-1}yw$ with $l\in \N$ and $w\in \fh^{1}$ varying under the condition $l+|w|=\textrm{weight}$.
Computations are performed by Mathematica. 
\vspace{2ex}
\begin{table}[!h]
\begin{center}
\caption{Number of Independent Derivation Relations for FMZVs}
\begin{tabular}{|c|r|r|r|r|r|r|r|r|r|r|r|r|r|} 
\hline
weight & 2& 3& 4& 5& 6& 7& 8& 9& 10& 11& 12& 13& 14 \\ 
\hline
relations & 1& 2& 5& 10& 22& 44& 90& 181& 363& 727& 1456& 2912& 5825 \\  
\hline
\end{tabular}
\end{center}
\end{table}

The interesting fact is that the number of independent relations of derivation relations in Table $1$ coincides with that of the original derivation relations in Table $2$, except that the weight is shifted by one. 
The reason for this coincidence is seen as follows. 
Write an element $w\in\fh^0$ as $w=xw', w'\in\fh^1$. 
Then by $\d_l(w)=xz^{l-1}yw'+x\d_l(w')$, the original derivation relations $\mathit{Z}(\d_l(w))=0$ can be written as 
\[ \mathit{Z}(x(\d_l(w')+z^{l-1}yw'))=0. \]
Hence the relation $\ZF (\d_l(w')+z^{l-1}yw')=0$ in weight $k$ exactly corresponds to the relation $\mathit{Z}(x(\d_l(w')+z^{l-1}yw'))=\mathit{Z}(\d_l(w))=0$ in weight $k+1$. 
\vspace{2ex}
\begin{table}[!h]
\begin{center}
\caption{Number of Independent Derivation Relations for MZVs}
\begin{tabular}{|c|r|r|r|r|r|r|r|r|r|r|r|r|r|} 
\hline
weight & 3& 4& 5& 6& 7& 8& 9& 10& 11& 12& 13& 14& 15 \\ 
\hline
relations & 1& 2& 5& 10& 22& 44& 90& 181& 363& 727& 1456& 2912& 5825 \\  
\hline
\end{tabular}
\end{center}
\end{table}

\section{Proofs of Theorem \ref{derivation0} and Theorem \ref{derivation1}}
We prove Theorem \ref{derivation0} by induction on $n=|w|$. 

\noindent (I) When $n=0$, i.e., $w=1$, we need to show 
\begin{align*}
&-\ZF (z^{l-1}y z^{m_{1}-1}y \cdots z^{m_{s}-1}y) \\
&+\sum_{i=1}^{s} \ZF (z^{m_1-1}y \cdots z^{m_{i-1}-1}y z^{m_i-1}x z^{l-1}y z^{m_{i+1}-1}y \cdots z^{m_{s}-1}y)= 0   
\end{align*}
for every $s\geq 0$. 
When $s=0$, by Theorem \ref{thm:duality}, we have
\[ -\ZF (z^{l-1}y)=\ZF (x^{l-1}y)=0. \]
Here, we note that $\zff (l)=0$ for any $l\in \N$. 
When $s\geq 1$, by Theorem \ref{3.3} and Theorem \ref{thm:duality}, 
\begin{align*}
&-\ZF (z^{l-1}y z^{m_{1}-1}y \cdots z^{m_{s}-1}y) \\
&+\sum_{i=1}^{s} \ZF (z^{m_1-1}y \cdots z^{m_{i-1}-1}y z^{m_i-1}x z^{l-1}y z^{m_{i+1}-1}y \cdots z^{m_{s}-1}y) \\
&=\ZF (-z^{l-1}yz^{m_1-1}y\cdots z^{m_s-1}y +z^{m_1-1}xz^{l-1}yz^{m_2-1}y\cdots z^{m_s-1}y \\
&\quad\quad\,\,\,\, + \quad\cdots\cdots\quad +z^{m_1-1}y\cdots z^{m_s-1}xz^{l-1}y) \\
&=(-1)^{s} \ZF (x^{l-1}yx^{m_1-1}y\cdots x^{m_s-1}y +x^{m_1-1}zx^{l-1}yx^{m_2-1}y\cdots x^{m_s-1}y \\
&\quad\quad\quad\quad\,\,\,\,\,\,\, +\quad\cdots\cdots\quad +x^{m_1-1}y\cdots x^{m_s-1}zx^{l-1}y) \\ 
&=(-1)^{s} \ZF (x^{m_1-1}y\cdots x^{m_s-1}y\ast x^{l-1}y) \\
&=(-1)^{s} \ZF (x^{m_1-1}y\cdots x^{m_s-1}y)\ZF (x^{l-1}y)=0. 
\end{align*}
%We also note that $\zff (k)=0$ for any $k\in \N$. 
 
\noindent (II) We assume the identity holds for $|w|=0,\ldots, n-1$ and for every $s\geq 0$. 
Suppose $w$ is of degree $n$. 
We may assume that $w$ is of the form $w=z^{r-1}yw'$ with $1\leq r\leq n, w'\in\fh^1$, by replacing $x^{r-1}y$ by $(z-y)^{r-1}y$ if $w$ starts with $x^{r-1}y$. 
\begin{align*}
\textrm{L.H.S.}&=Z_{\mathcal{F}} (z^{m_1-1}y\cdots z^{m_s-1}y\partial_{l}(z^{r-1}yw')) \\
&=Z_{\mathcal{F}} (-z^{m_1-1}y\cdots z^{m_s-1}yz^{r-1}xz^{l-1}yw' +z^{m_1-1}y\cdots z^{m_s-1}yz^{r-1}y\partial_{l}(w')). 
\end{align*}
By the induction hypothesis, we have
\begin{align*}
Z_{\mathcal{F}} (z^{m_1-1}y\cdots z^{m_s-1}yz^{r-1}y\partial_{l}(w')) &=Z_{\mathcal{F}} (-z^{l-1}yz^{m_1-1}y\cdots z^{m_s-1}yz^{r-1}yw' \\
&\quad\quad\,\,\,\,\, +z^{m_1-1}xz^{l-1}yz^{m_2-1}y\cdots z^{m_s-1}yz^{r-1}yw' \\ 
&\quad\quad\,\,\,\,\, +\quad\cdots\cdots \\
&\quad\quad\,\,\,\,\, +z^{m_1-1}y\cdots z^{m_s-1}yz^{r-1}xz^{l-1}yw'). 
\end{align*}
Thus, 
\begin{align*}
\textrm{L.H.S.}&=\ZF (-z^{m_1-1}y\cdots z^{m_s-1}yz^{r-1}xz^{l-1}yw' -z^{l-1}yz^{m_1-1}y\cdots z^{m_s-1}yz^{r-1}yw' \\
&\quad\quad\,\,\,\,\, +z^{m_1-1}xz^{l-1}yz^{m_2-1}y\cdots z^{m_s-1}yz^{r-1}yw' +\quad\cdots\cdots \\
&\quad\quad\,\,\,\,\, +z^{m_1-1}y\cdots z^{m_s-1}xz^{l-1}yz^{r-1}yw' +z^{m_1-1}y\cdots z^{m_s-1}yz^{r-1}xz^{l-1}yw') \\
&=\ZF (-z^{l-1}yz^{m_1-1}y\cdots z^{m_s-1}yz^{r-1}yw' +z^{m_1-1}xz^{l-1}yz^{m_2-1}y\cdots z^{m_s-1}yz^{r-1}yw' \\ 
&\quad\quad\,\,\,\,\, +\quad\cdots\cdots\quad +z^{m_1-1}y\cdots z^{m_s-1}xz^{l-1}yz^{r-1}yw') \\
&=\textrm{R.H.S.}, 
\end{align*}
and hence the identity holds for $n$ and by induction, the proof is done.  

\bigskip

Now, we prove Theorem \ref{derivation1} by induction on $t$. 
We have proved the case $t=1$. 
We assume the identity holds when the number of derivations on the left is less than $t$. 
\begin{align*}
&\ZF (z^{m_1-1}y\cdots z^{m_s-1}y\d_{l_1}\cdots\d_{l_t}(w)) \\ 
&=\ZF (-z^{l_1-1}yz^{m_1-1}y\cdots z^{m_s-1}y\d_{l_2}\cdots\d_{l_t}(w) \\
&\quad\quad\,\,\,\,\, +z^{m_1-1}xz^{l_1-1}yz^{m_2-1}y\cdots z^{m_s-1}y\d_{l_2}\cdots\d_{l_t}(w) \\ 
&\quad\quad\,\,\,\,\, +\quad\cdots\cdots\quad \\
&\quad\quad\,\,\,\,\, +z^{m_1-1}yz\cdots z^{m_s-1}xz^{l_1-1}y\d_{l_2}\cdots\d_{l_t}(w)) \\ 
&=\ZF (-z^{l_1-1}yz^{m_1-1}y\cdots z^{m_s-1}y\d_{l_2}\cdots\d_{l_t}(w) \\
&\quad\quad\,\,\,\,\, +z^{m_1-1}(z-y)z^{l_1-1}yz^{m_2-1}y\cdots z^{m_s-1}y\d_{l_2}\cdots\d_{l_t}(w) \\ 
&\quad\quad\,\,\,\,\, +\quad\cdots\cdots\quad \\
&\quad\quad\,\,\,\,\, +z^{m_1-1}yz\cdots z^{m_s-1}(z-y)z^{l_1-1}y\d_{l_2}\cdots\d_{l_t}(w)) \\ 
&=(-1)^{s} \sum_{i=0}^{s} \sum_{\sigma\in S_{s+t}^{(s+1)}} \ZF (z^{a'_{i,\sigma (1)}-1}u^{\sigma}_1 \cdots z^{a'_{i,\sigma (s+t)}-1} u^{\sigma}_{s+t}w) \\
&\quad +(-1)^{s} \sum_{i=1}^{s} \sum_{\sigma\in S_{s+t-1}^{(s)}} \ZF (z^{a''_{i,\sigma (1)}-1}u^{\sigma}_1 \cdots z^{a''_{i,\sigma (s+t-1)}-1} u^{\sigma}_{s+t-1}w), 
\end{align*}
where $\textrm{\textup{\textmd{\textbf{a}}}}'_i=(a'_{i,1},\ldots ,a'_{i,s+t})=(m_1,\ldots ,m_i,l_1,m_{i+1},\ldots ,m_s,l_2,\ldots ,l_t)$ and $\textrm{\textup{\textmd{\textbf{a}}}}''_i=(a''_{i,1},\ldots ,a''_{i,s+t-1})=(m_1,\ldots ,m_{i-1},m_i+l_1,m_{i+1},\ldots ,m_s,l_2,\ldots ,l_t)$ in the last summation. 
We let 
\begin{align*}
L:&=\sum_{i=0}^{s} \sum_{\sigma\in S_{s+t}^{(s+1)}} \ZF (z^{a'_{i,\sigma (1)}-1}u^{\sigma}_1 \cdots z^{a'_{i,\sigma (s+t)}-1} u^{\sigma}_{s+t}w), \\
M:&=\sum_{i=1}^{s} \sum_{\sigma\in S_{s+t-1}^{(s)}} \ZF (z^{a''_{i,\sigma (1)}-1}u^{\sigma}_1 \cdots z^{a''_{i,\sigma (s+t-1)}-1} u^{\sigma}_{s+t-1}w), \\ 
N:&=\sum_{\sigma\in S_{s+t}^{(s)}} \ZF (z^{a_{\sigma (1)}-1}u^{\sigma}_1 \cdots z^{a_{\sigma (s+t)}-1} u^{\sigma}_{s+t}w). 
\end{align*} 
For each element in $L$, there exists a unique element in $N$ such that they are corresponding to each other except for one letter $u_i$ between $z^{m_i-1}$ and $z^{l_1-1}$, which is $-y$ in $L$ and $x$ in $N$. 
Similarly, by understanding $z^{m_i+l_1-1}=z^{m_i-1}\cdot z\cdot z^{l_1-1}$, there is one-to-one correspondence between the elements in $M$ and $N$ such that they are corresponding to each other except for $u_i$ between $z^{m_i-1}$ and $z^{l_1-1}$, which is $z$ in $M$ and $x$ in $N$. 
Since $x=-y+z$, we have 
\begin{align*} 
& \ZF (z^{m_1-1}y \cdots z^{m_s-1}y \partial_{l_1} \cdots \partial_{l_t}(w)) \\
& = (-1)^{s} \sum_{\sigma\in S_{s+t}^{(s)}} \ZF (z^{a_{\sigma (1)}-1}u^{\sigma}_1 \cdots z^{a_{\sigma (s+t)}-1} u^{\sigma}_{s+t}w) \quad (w \in \fh^1). 
\end{align*} 
Then, we find the identity holds for $t$. 

\section*{Acknowledgment}
The author would like to express his sincere gratitude to Professor Masanobu Kaneko for valuable comments and advices. 
He would also like to thank Kojiro Oyama for interesting suggestions that inspired this work.

\end{document}